# DOMAIN OF ATTRACTION OF ASYMPTOTICALLY STABLE PERIODIC SOLUTIONS OBTAINED VIA AVERAGING PRINCIPLE

Oleg Makarenkov

*Abstract:* In this paper we propose an approach to evaluate the domain of attraction of asymptotically stable periodic solutions obtained via averaging principle (second Bogolubov's theorem or Mel'nikov's method). We discuss also how this result is extended in the case when the right hand part is nonsmooth.

**1. Introduction**

Since the classical paper by van der Pol [6] many oscillating process in mechanics and engineering are described by means of the following system of ordinary differential equations

$$\dot{u} = A(u) + \varepsilon h(t,u,\varepsilon), \qquad u \in \mathrm{R}^n, \tag{1}$$

where $\varepsilon > 0$ is a small parameter, $A \in C^2(\mathrm{R}^n,\mathrm{R}^n)$, $h \in C^1(\mathrm{R}\times\mathrm{R}^n\times[0,1],\mathrm{R}^n)$ is a $T$-periodic in time and any solution of the unperturbed system $\dot{u} = A(u)$ is $T$-periodic. A well known change of variables allows us to rewrite (1) in the standard form of Krylov-Bogolubov averaging principle (see [5], §2):

$$\dot{x} = \varepsilon g(t,x,\varepsilon), \qquad x \in \mathrm{R}^n, \tag{2}$$

where $g \in C^1(\mathrm{R}\times\mathrm{R}^n\times[0,1],\mathrm{R}^n)$ and analysis of $T$-periodic oscillations in (1) is then reduced to the same analysis for (2). A key role in this analysis is played by the averaging function:

$$g_0(v) = \int_0^T g(\tau,v,0)d\tau,$$

namely the Second Bogolubov's theorem claims:

**Second Bogolubov's theorem** ([2], Ch.1, §5, Theorem II). *Properties $g_0(v_0) = 0$ and $\det(g_0)'(v_0) \neq 0$ assure the existence and uniqueness, for $\varepsilon > 0$ small of a $T$-periodic solution $x_\varepsilon$ of system* (2) *in a neighborhood of $v_0$, while the fact that all the eigenvalues of $(g_0)'(v_0)$ have negative real part, provides also its asymptotic stability.*

The goal of the present paper is to answer the following question settled in a similar situation by Andronov-Vitt-Khaikin ([1], Ch.III, §3): *What is the domain of attraction of the asymptotically stable T-periodic solution $x_\varepsilon$ given by the Second Bogolubov's theorem?* Surprisingly, we were not able to find the answer in the literature.

## 2. Main result

We can give the following answer to the question posted in the introduction.

**Theorem 1.** *Let $g \in C^1(R \times R^n \times [0,1], R^n)$ and $v_0 \in R^n$ satisfies $g_0(v_0)=0$. Assume, that the convex set $v_0 \in V \subset R^n$ is chosen in such a way that there exist constants $\alpha > 0$, $q \in [0,1]$ and a norm $||\cdot||_0$ on $R^n$ satisfying the property*

$$\|v_1 + \alpha g_0(v_1) - v_2 - \alpha g_0(v_2)\|_0 \leq q \|v_1 - v_2\|_0 \quad \text{for any} \quad v_1, v_2 \in V. \tag{3}$$

*Then for any $\varepsilon > 0$ sufficiently small system (2) has an unique asymptotically stable T-periodic solution $x_\varepsilon$ with $x_\varepsilon(0) \in V$. Moreover, $x_\varepsilon$ attracts any other solution with initial condition in V.*

**Proof.** Denote by $x(\cdot, v, \varepsilon)$ the *T*-periodic solution of (2) with initial condition $x(0)=v$. We have

$$x(t,v,\varepsilon) = v + \varepsilon \int_0^T g(s, x(s,v,\varepsilon), \varepsilon) d\tau = v + \varepsilon g_0(v) + \varepsilon(g_\varepsilon(v) - g_0(v)), \tag{4}$$

where $g_0(v) = \int_0^T g(\tau, v, 0) d\tau$ and $g_\varepsilon(v) = \int_0^T g\left(\tau, v + \varepsilon \int_0^\tau g(s, x(s,v,\varepsilon), \varepsilon) ds, \varepsilon \right) d\tau$.

If can be easily checked (see e.g. [3], Lemma 2.3) that (3) implies

$$\|v_1 + \varepsilon g_0(v_1) - v_2 - \varepsilon g_0(v_2)\|_0 \leq (1 - \varepsilon(1-q)\alpha) \|v_1 - v_2\|_0, \quad \text{for any} \quad v_1, v_2 \in V, \varepsilon \in (0, \varepsilon_0]. \tag{5}$$

Choose $r > 0$ sufficiently small such that $v_0 \in B_r(v_0) \subset V$. Then from (5) we have that there exists $\alpha > 0$ such that $I + \alpha g_0$ maps $B_r(v_0)$ into itself. Therefore, we can decrease $\varepsilon_0 > 0$, if necessary, in such a way that $I + \alpha g_0$ maps $B_r(v_0)$ into itself for any $\varepsilon \in (0, \varepsilon_0]$ as well. By the Brouwer theorem (see, e.g. [4], Theorem 3.1) we have that $B_r(v_0)$ contains at least one fixed point of the map $I + \alpha g_0$ for any $\varepsilon \in (0, \varepsilon_0]$. Denote this fixed point by $v_\varepsilon$. Than we have $g_\varepsilon(v_\varepsilon) = 0$ and $x(T, v_\varepsilon, \varepsilon) = v_\varepsilon$ for any $\varepsilon \in (0, \varepsilon_0]$.

Now we prove that $x(\cdot, v_\varepsilon, \varepsilon)$ is asymptotically stable *T*-periodic solution of (2) globally in *V*. Since $(v, \varepsilon) \to g_\varepsilon(v)$ is continuously differentiable then $\varepsilon_0 > 0$ can be diminished, if necessary, in such away that

$$\|g_\varepsilon(v_1) - g_0(v_1) - g_\varepsilon(v_2) + g_0(v_2)\|_0 \leq (1-q)(\alpha/2)\|v_1 - v_2\|_0, \quad \text{for any} \quad v_1, v_2 \in V, \varepsilon \in (0, \varepsilon_0]. \tag{6}$$

Combining (4), (5) and (6) we have

$$\|x(T,v_1,\varepsilon) - x(T,v_2,\varepsilon)\|_0 \leq (1-\varepsilon(1-q)\alpha +\varepsilon(1-q)(\alpha/2)) \|v_1 - v_2\|_0 \quad \text{for any} \quad v_1, v_2 \in V, \varepsilon \in (0,\varepsilon_0].$$

Since $0 < 1-\varepsilon(1-q)\alpha +\varepsilon(1-q)(\alpha/2) < 1$ then the latter implies (see [4], Lemma 9.2) the thesis of the lemma.

□

**Remark 1.** *Assumptions of the Second Bogolubov's theorem imply the existence of $\alpha > 0$, $q \in [0,1]$ and a norm $\|\cdot\|_0$ on $R^n$ satisfying (3) of Theorem 1 with sufficiently small ball centered at $v_0$ playing the role of V (see [3], Lemma 2.7).*

**3. An extension to nonsmooth differential equations**

Following the lines of [3] it is possible to give a generalization of Theorem 1 for the case when the right hand part in (2) is only Lipschitzian, but possesses the following property:

**Definition.** *Non-autonomous function $(t,v,\varepsilon) \to g(t,v,\varepsilon)$ (with time variable t) depending on parameter $\varepsilon \in [0,1]$ is said to be weakly differentiable if given any $\gamma > 0$ there exist $\delta > 0$ and a family*

$$\left\{ \bigcup_{i=1}^{N(v)} [s_{i,v}, t_{i,v}] \right\}_{v \in V} \subset [0,T] \text{ satisfying } s_1 \leq t_1 \leq s_2 \leq t_2 \leq \ldots, \sum_{i=1}^{N(v)} |s_{i,v} - t_{i,v}| \leq \gamma \text{ such that for any } v \in V,$$

$t \in [0,T] \setminus \bigcup_{i=1}^{N(v)} [s_{i,v}, t_{i,v}]$ *and $\varepsilon \in [0,1]$ we have*

   a) *the function $g(t,\cdot,\varepsilon)$ is differentiable in $\delta$-neighborhood of v,*

   b) $\| g'_v(t, v + h_1, \varepsilon) - g'_v(t, v + h_2, \varepsilon) \| \leq \gamma$ *for any* $\|h_1\| \leq \delta, \|h_2\| \leq \delta$.

One can prove the following result.

**Theorem 2.** *Let $g \in C^0(R \times R^n \times [0,1], R^n)$ be locally Lipschitzian in the second variable weakly differentiable function and $v_0 \in R^n$ satisfies $g_0(v_0)=0$. Assume, that the convex set $v_0 \in V \subset R^n$ is chosen in such a way that there exist constants $\alpha > 0$, $q \in [0,1]$ and a norm $\|\cdot\|_0$ on $R^n$ satisfying the property*

$$\|v_1+\alpha g_0(v_1)-v_2-\alpha g_0(v_2)\|_0 \leq q\|v_1-v_2\|_0 \quad \text{for any} \quad v_1, v_2 \in V. \tag{3}$$

*Then for any $\varepsilon > 0$ sufficiently small system (2) has an unique asymptotically stable T-periodic solution $x_\varepsilon$ with $x_\varepsilon(0) \in V$. Moreover, $x_\varepsilon$ attracts any other solution with initial condition in V.*

**Remark 2.** *By several particular examples we convinced ourselves that if function $(t,v,\varepsilon) \to h(t,v,\varepsilon)$ in (1) is piece-wise differentiable in the second variable then $(t,v,\varepsilon) \to g(t,v,\varepsilon)$ in (2) is weekly differentiable in the sense of above definition, but we did not prove a general statement up to now.*

**Acknowledgements.** The work is partially supported by the grant the Grant BF6M10 of Russian Federation Ministry of Education and CRDF (US), and by RFBR Grants 06-01-72552, 05-01-00100.

**References**

[1] A. A. Andronov, A. A. Vitt, S. E. Khaikin, Theory of oscillators. Translated from the Russian by F. Immirzi; translation edited and abridged by W. Fishwick Pergamon Press, Oxford-New York-Toronto, Ont. 1966, 815 pp.

[2] N. N. Bogolyubov, On Some Statistical Methods in Mathematical Physics, Akademiya Nauk Ukrainskoi SSR, 1945, 139 pp. (Russian)

[3] A. Buica, J. Llibre, O. Makarenkov, Asymptotic stability of periodic solutions for nonsmooth differential equations with application to the nonsmooth van der Pol oscillator, SIAM J. Math. Anal., submitted.

[4] M. A. Krasnosel'skii, The operator of translation along the trajectories of differential equations, Translations of Mathematical Monographs, 19. Translated from the Russian by Scripta Technica, American Mathematical Society, Providence, R.I. 1968.

[5] N. Kryloff , N. Bogoliouboff, La theorie generale de la mesure dans son application a l'etude des systemes dynamiques de la mecanique non lineaire. (French) Ann. of Math. (2) 38 (1937), no. 1, 65–113.

[6] B. Van der Pol, On relaxation oscillations, Philosophical Magazine Series 7, 2 (1926) 978–992.

Oleg Makarenkov,
Research Institute of Mathematics, Voronezh State University, Russia,
394006, Russia, Voronezh, Universitetskaja pl. 1,
omakarenkov@math.vsu.ru